\documentclass[final]{siamart1116}

\usepackage[utf8]{inputenc}
\usepackage[T1]{fontenc}
\usepackage{amsfonts}
\usepackage{grffile}
\usepackage{tikz}
\usepackage{pgfplots}
\usepackage{pgfplotstable}
\pgfplotsset{compat=1.10}
\usepgfplotslibrary{fillbetween}
\usepackage{eurosym}
\usepackage{graphicx}
\usepackage{color}
\usepackage{listings}
\usepackage{paralist}
\usepackage{curves}
\usepackage{calc}
\usepackage{picinpar}
\usepackage{enumerate}
\usepackage{algorithm}
\usepackage{algpseudocode}
\usepackage{bm}
\usepackage{multibib}
\usepackage{hyperref}
\usepackage{textcase}
\usepackage{nicefrac}
\usepackage{stmaryrd}
\usepackage{tabu,booktabs}
\usepackage{todonotes}
\usepackage{tabularx}
\usepackage{dcolumn}

\usepackage[detect-weight,binary-units]{siunitx}

\DeclareSIUnit\dof{DOF}
\DeclareSIUnit\flop{FLOP}

\newcolumntype{Y}{%
  S[
  table-format=1.2e1,
  table-auto-round,
  scientific-notation=true,
  ]}

\newcolumntype{F}{%
  S[
  table-format=2.1,
  round-mode=figures,
  round-precision=3,
  ]}

\newcolumntype{G}{%
  S[
  table-format=3,
  round-mode=figures,
  round-precision=3,
  ]}

\RenewDocumentCommand{\vec}{m}
{\bm{{#1}}}

\definecolor{listingbg}{gray}{0.95}

\theoremstyle{definition}

\numberwithin{theorem}{section}

\colorlet{texcscolor}{blue!50!black}
\colorlet{texemcolor}{red!70!black}
\colorlet{texpreamble}{red!70!black}
\colorlet{codebackground}{black!25!white!25}
\lstdefinestyle{siamlatex}{%
  style=tcblatex,
  texcsstyle=*\color{texcscolor},
  texcsstyle=[2]\color{texemcolor},
  keywordstyle=[2]\color{texemcolor},
  moretexcs={cref,Cref,maketitle,mathcal,text,headers,email,url},
}

\title{High-performance Implementation of Matrix-free High-order Discontinuous Galerkin Methods%
\thanks{Submitted to the editors 2017--11-29.
\funding{This project was supported in part by Deutsche Forschungsgemeinschaft under grant Ba 1498/10-2 within
the SPPEXA programme.
For computational resources the authors acknowledge support by the state of Baden-Württemberg through bwHPC
and the German Research Foundation (DFG) through grant INST 35/1134-1 FUGG.}}}
\author{Steffen Müthing\thanks{Heidelberg University, Interdisciplinary Center for Scientific Computing,
Im Neuenheimer Feld 205, D-69120 Heidelberg, Germany. \newline
\texttt{\{steffen.muething,marian.piatkowski,peter.bastian\}@iwr.uni-heidelberg.de}
} \and Marian Piatkowski$^\dagger$ \and Peter Bastian$^\dagger$}
\date{\today}

\begin{document}

\maketitle

\begin{abstract}
Achieving a substantial part of peak performance on todays and future high-performance computing
systems is a major challenge for simulation codes. In this paper we address this question in the context of
the numerical solution of partial differential equations with finite element methods, in particular
the discontinuous Galerkin method applied to a convection-diffusion-reaction model problem.
Assuming tensor product structure of basis functions and quadrature on cuboid meshes in
a matrix-free approach a substantial
reduction in computational complexity can be achieved for operator application compared to a
matrix-based implementation while at the same time enabling SIMD vectorization and the use of
fused-multiply-add.  Close to 60\% of peak performance
are obtained for a full operator evaluation on a Xeon Haswell CPU with 16 cores and speedups
of several hundred (with respect to matrix-based computation)
are achieved for polynomial degree seven.
Excellent weak scalability on a single node as well as the roofline model demonstrate that the
algorithm is fully compute-bound with a high flop per byte ratio. Excellent scalability is also
demonstrated on up to 6144 cores using message passing.

%

Keywords: High-order discontinuous Galerkin methods $\cdot$ matrix-free methods $\cdot$ sum factorization
$\cdot$ SIMD vectorization
\end{abstract}

\section{Introduction}

Achieving a substantial part of peak performance on todays and future high-performance computing
systems is a major challenge for simulation codes. In this paper we address this question in the context of
the numerical solution of partial differential equations (PDEs) with grid-based methods, in particular finite element
methods. This field comprises a good part of todays supercomputer applications.
As many supercomputers are based on standard components such as multi-core CPUs
the techniques discussed in this paper are also relevant on desktop systems.

In the early 2000s the almost effortless increase in performance due to increasing
clock rate ended rather abruptly \cite{Asanovic:2009:VPC:1562764.1562783}. As a result, the
major challenges for simulation software arising from computer architecture
today are \cite{Dongarra01022011,Keyes:2011:ETW,Dongarra:2014:AMR}:
\begin{itemize}
\item \textit{Massive Parallelism.} The increase in transistor count per chip has been used to place multiple independent cores
(typically 4 to 64) on a CPU chip.
Peak performance of such a core is only obtained through the use of SIMD (vector-)
instructions where one instruction (such as a fused-multiply-add)
is executed on several (say 4 to 16) operands in parallel.
In a supercomputer a fixed small number
of such CPUs with shared memory are combined in a node and many such nodes are connected via a scalable message passing network.
An example of such a system is the 100 Petaflop system TaihuLight \cite{TaihuLight2016}.
As a consequence, applications on such systems therefore need to harness massive parallelism on the order of $10^6$ to
$10^9$ floating point operations executed simultaneously
on three different levels accessible via three different programming models: (i) explicit
SIMD instructions, (ii) multiple threads using shared memory and (iii) parallel processes communicating via message passing.
\item \textit{Memory gap.} Main (DRAM) memory is not able to deliver operands at the speed as they are processed
in the CPU. This fact has always existed in computer architecture \cite{Wilkes:2001:MGF:373574.373576} and the introduction
of multicore CPUs has brought no relief as the number of memory controllers is typically (much) smaller
than the number of cores. Moreover, data transfer from main memory comprises the main energy consumption of a computation.
In \cite{Keyes:2011:ETW} it was estimated that only 0.1 byte per floating point operation (flop) can be afforded to stay
within the power envelope of a future exascale system (about 20 MW). This corresponds to 80 flops
to be executed per double precision floating point number loaded from memory. As a consequence
it is of great concern to develop algorithms and implementations maximizing the amount of
\textit{useful} flops per byte transfered from/to main memory.
\item \textit{Heterogeneity.} Many current supercomputing systems employ heterogeneous architectures combining
multicore CPUs with one or several coprocessors such as a GPU or the Intel MIC. Coprocessors implement an
increased number of simplified cores combined with higher bandwidth to a dedicated memory. Programming such systems
requires careful assignment of code parts to CPU and coprocessor and managing the necessary data transfer.
Software layers such as OpenACC, OpenCL and Kokkos aim at perfomance-portable programming of heterogeneous systems
which is complicated by the fact that different systems often need different data layout.
In this paper we focus on high-performance CPU implementations and do not consider coprocessor architectures
explicitely. However, the methodology developed is, in principle, transferable to coprocessor architectures
and therefore is also relevant on these architectures.
\end{itemize}

Finite element (FE) discretization of stationary, nonlinear PDEs (or implicit in time discretization of
instationary PDEs) leads
to the solution of large and sparse nonlinear algebraic systems of the form $R(z)=0$.
Iterative solution then leads, after linearization
(or if the system is linear to begin with), to large sparse linear systems $Az=b$, see e.g. \cite{Ern,Elman2005}.
A basic ingredient of iterative solvers is residual evaluation
$r^k=R(z^k)$ and, respectively, operator application $y^k=A z^k$.
Note that explicit in time discretization leads to a basic operation of the form
$z^{k+1} = z^{k} + F(z^k)$ where $F(z^k)=M^{-1}R(z^k)$ with $M$ the mass matrix.
In particular, when the
mass matrix is diagonal this operation can be viewed as a variant of residual evaluation.
A single matrix-vector product (with stored matrix $A$) is
memory bound as only two flops are executed per matrix element read from memory
(this holds for sparse and dense matrices). Additional complications in the case of sparse matrices such as index access and indirect memory access
as well as limited cache reuse on the vectors leads to poor floating point performance (relative to peak performance)
of stored matrix vector products even in highly optimized formats \cite{Liu:2013:ESM:2464996.2465013,doi:10.1137/130930352}.
The difficulty of achieving a substantial part of peak performance with implicit, finite-element based PDE solvers
is illustrated by the 2015 Gordon Bell prize winner providing the state-of-the-art.
The authors of \cite{Rudi:2015:EIS:2807591.2807675} state that their
implicit multigrid solver applied to a highly nonlinear earth mantle convection problem is memory bound
(despite using a matrix-free high-order discretization)
and report in figure 7 an average performance of
7.5 GFlops/s per node which corresponds to 3.5\% of peak performance of a BlueGene/Q node (204,5 GFlops/s).

Matrix-free operator application (or residual evaluation and linearized operator
application in the nonlinear case) is a promising technique to
increase the flop per byte ratio substantially and thus to overcome the memory gap.
Obviously, to be faster in real time, this
requires to be able to compute matrix entries faster than to load them from memory. For low order
discretizations, this approach is popular with stencil-based codes and constant coefficients
\cite{Douglas2000,Sellappa:2004,Datta:2008} but has
also been extended successfully to certain types of unstructured meshes e.g.~in the HHG software \cite{doi:10.1137/130941353}.
The situation is fundamentally different for high-order (spectral) finite element methods. Here, the
sum-factorization technique for a tensor product basis leads to a \textit{significant reduction of computational
complexity} of matrix-free operator evaluation compared to a matrix-vector product
\cite{ORSZAG198070,Karniadakis2005,Kronbichler2012135}.
Moreover, the operations can be arranged in a sequence of matrix-matrix products
\cite{Buis:1996:EVP:225545.225548} which lend themselves for vectorized execution of
fused-multiply-add (fma) instructions.
Thus, less operations are executed
at a higher rate \cite{Kronbichler2012135}.
Sum-factorization is most easily applied to cuboid elements but has also been used
for simplicial elements \cite{Karniadakis2005,Joachim2014}. It can handle nonlinear problems as well as
higher-order geometries and has been used to set up Jacobians \cite{Melenk20014339}.
Several established simulation software packages for solving PDEs implement high-order
finite element methods using sum factorization, among them
nek5000\footnote{\url{http://www.mcs.anl.gov/project/nek5000-computational-fluid-dynamics-code}},
Nektar++\footnote{\url{https://www.nektar.info}},
deal.II\footnote{\url{https://www.dealii.org}} and
NGSolve\footnote{\url{https://ngsolve.org}}. Here we report on our implementation within
the DUNE\footnote{\url{https://www.dune-project.org}} software framework \cite{Dune2008a,Dune2008b}.

\begin{figure}
\begin{center}
\begin{tikzpicture}[scale=0.85]
\begin{axis}[
xlabel={Polynomial Degree},
ylabel={FLOPs per DOF},
grid=major,
legend entries={{volume, $d=2$},{volume, $d=3$},{face, $d=2$},{face, $d=3$}},
legend pos=outer north east,
]
\addplot[red,domain=1:8,samples=8,mark=*,]{(4*3*2*(x+1)+2*2*2+5*2+3)};
\addplot[green,domain=1:8,samples=8,mark=+,]{(4*4*3*(x+1)+2*3*3+5*3+3)};
\addplot[blue,domain=1:8,samples=8,mark=*,]{8*3*2*2+2*(2*2*2+5*2+3)/(x+1)};
\addplot[magenta,domain=1:8,samples=8,mark=+,]{8*4*3*3+2*(2*3*3+5*3+3)/(x+1)};
\end{axis}
\end{tikzpicture}
\end{center}
\caption{Number of floating point operations per degree of freedom for sum-factorized volume and face
terms in discontinuous Galerkin on a structured mesh in two and three space dimensions.}
\label{fig:flops_per_dof}
\end{figure}
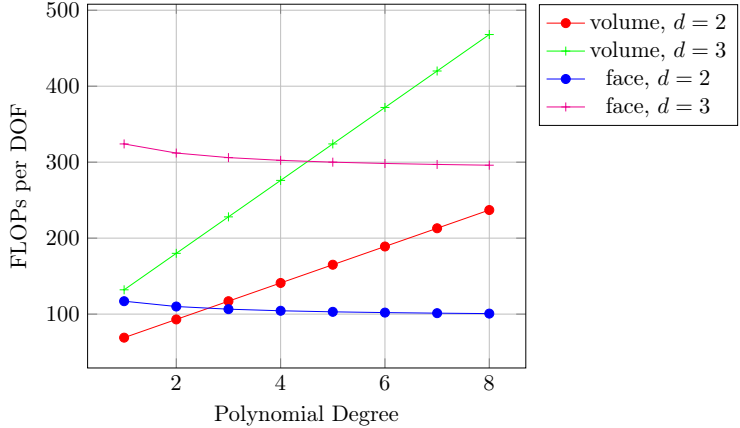

In this paper we combine sum-factorization with discontinuous Galerkin (DG) finite element methods.
DG methods are popular in the CFD community due to their local mass conservation property while
being higher order convergent (for sufficiently regular problems) and their ability to handle
elliptic, parabolic and hyperbolic problems \cite{DGProceedings00}. In particular we are interested in
porous media flow applications \cite{Ern20101491,Huang2013169,ba_amg_twophase_2014,Li2015107}
but the methodology developed here is even more generally applicable. When compared to continuous Galerkin (CG)
methods DG involves more degrees of freedom for the same mesh and polynomial degree. This could be alleviated
with hybrid DG methods \cite{doi:10.1137/070706616} or hybrid high-order methods
\cite{DIPIETRO201531}, but we restrict ourselves to the symmetric
interior penalty DG method with weighted averages (WSIPG) \cite{Ern01042009} which performs well for elliptic problems
with highly varying, anisotropic diffusion coefficients \cite{fvca6_2011}. A significant advantage of DG
methods over CG methods as it turns out is that degrees of freedom can be stored consecutively per element and
no gather/scatter operations are necessary to evaluate the element and face integrals.
On the other hand, the challenge in DG methods are the interior face integrals where sum-factorization
is less efficient due to the lower dimension and the flop per byte ratio is worse. Figure \ref{fig:flops_per_dof} illustrates that
computational work for evaluating face integrals is dominating over the work in the volume integrals up to
polynomial degree three in 2$d$ and five in 3$d$. Moreover,
more data has to be moved per flop in the face terms.

In addition to operator evaluation, for many problems efficient, robust and scalable preconditioners
are essential to achieve fast convergence of iterative methods. While being essential, matrix-free preconditioners
are beyond the scope of this paper. We will treat them in a forthcoming paper based on our previous
work on low-order subspace correction using algebraic multigrid (AMG) \cite{amg4dg}. All operations involving the high-order
DG system can be implemented matrix-free while AMG on the low-order subspace provides robustness with respect
to anisotropy and high coefficient contrasts. Others have shown good scalability and performance for
matrix-free geometric multigrid \cite{doi:10.1137/130941353} as well as matrix-free geometric multigrid combined with algebraic multigrid
\cite{Rudi:2015:EIS:2807591.2807675}. For this paper we have deliberately chosen to solely  concentrate
on matrix-free operator evaluation and matrix-free explicit time-stepping schemes in the context of
linear, variable coefficient PDEs.

The rest of the paper is organized as follows: In Section \ref{sec:problem_discretization}
we review the model PDE and its DG discretization, Section \ref{Sec:SpectralDG} describes the
sum-factorization technique for discontinuous Galerkin methods while Section \ref{Sec:Implementation}
provides implementation details before Section \ref{Sec:Results} presents
extensive numerical results and performance evaluation. Then we conclude in Section \ref{Sec:Conclusion}.

\section{Continuous Problem and Discontinuous Galerkin Formulation}
\label{sec:problem_discretization}

\subsection{Model Problem}

In this work we consider the convection-diffusion-reaction model problem in a finite domain
$\Omega\subset\mathbb{R}^d$ and time interval $\Sigma=(t_0,t_0+T_F]$:
\begin{subequations}\label{eq:modelproblem}
\begin{align}
\partial_t u + \nabla\cdot(\vec{b}u-D\nabla u) + cu &= f &&\text{in $\Omega\times\Sigma$}, \label{eq:problem}\\
u &= g &&\text{on $\Gamma_D\times\Sigma\subseteq\partial\Omega\times\Sigma$}, \label{eq:dirichletbc}\\
(\vec{b}u - D\nabla u)\cdot\vec{\nu} &= j &&\text{on $\Gamma_N\times\Sigma\subseteq\partial\Omega\times\Sigma$},\\
(D\nabla u)\cdot\vec{\nu} &= 0 &&\text{on $\Gamma_O\times\Sigma=(\partial\Omega\setminus
\Gamma_D\setminus\Gamma_N)\times\Sigma$},\label{eq:outflow}\\
u &= u_0 &&\text{at $t=t_0$}, \label{eq:initial_condition}
\end{align}
\end{subequations}
for the unknown function $u:\Omega\times\Sigma\to\mathbb{R}$.
The scalar, vector and matrix-valued coefficients
$c(\vec{x},t)$, $\vec{b}(\vec{x},t)$ and $D(\vec{x},t)$ are understood as functions of a spatial variable $\vec{x}\in\Omega$
and time $t\in\Sigma$, $\vec{\nu}(\vec{x})$ is the unit outer normal vector.
$D(\vec{x},t)$ is symmetric, uniformly positive definite and may be discontinuous
(with discontinuities resolved by the computational mesh). The velocity field $\vec{b}(\vec{x},t)$ is assumed to
be in $H(\text{div};\Omega)$ for fixed $t$.
The outflow boundary condition \eqref{eq:outflow} is only used
in the convection-dominated case and it is assumed that
$\Gamma_O\subseteq\Gamma^+ = \{\vec{x}\in\partial\Omega:
\vec{b}(\vec{x})\cdot\vec{\nu}(\vec{x})\geq 0\}$.
Below we will also consider the stationary variant of equation \eqref{eq:modelproblem}
where $\partial_t u = 0$, all coefficients depend only on space and the initial
condition \eqref{eq:initial_condition} is omitted.
This model problem occurs (at least as a component) in a wide range of applications
such as heat and fluid flow or computational biology.
We are particularly interested in porous media flow where (variants of)  equation
\eqref{eq:modelproblem} model the flow of one or several fluid phases
possibly coupled with transport of dissolved substances \cite{ba_amg_twophase_2014,Huang2013169,Li2015107}.
One difficulty of flows in porous media is that coefficients $\vec{b}$ and $D$ are highly varying in space and time.

\subsection{Discontinuous Galerkin Discretization}

As spatial discretization of equation \eqref{eq:modelproblem} we employ the
weighted symmetric interior penalty discontinuous Galerkin (WSIPG) method introduced in
\cite{Ern20101491}. This method is attractive as it is locally mass conservative,
works on very general meshes, allows varying polynomial degree,
can handle highly varying, anisotropic and
matrix-valued diffusion coefficients \cite{fvca6_2011} and is also
suited for the convection-dominated case \cite{Ngo2015331}.
$H(\text{div};\Omega)$ flow fields, if required, can be reconstructed locally
by means of interpolation \cite{ErnHdiv2007,BastianRiviere}.

Let $\{\mathcal{T}_h\}_{h>0}$ be a family of shape regular triangulations of the
domain $\Omega$ consisting of closed elements $T$, each being
the image of a map $\mu_T : \hat T \to T$ with $\hat T$ the reference
cube in $d$ dimensions. The map $\mu_T$ is differentiable, invertible and
its gradient is nonsingular on $\hat T$.
The diameter of $T$ is $h_T$ and $\vec{\nu}_T$ is its unit outer normal vector.
$F$ is an interior face if it is the intersection of two elements
$T^-(F), T^+(F)\in\mathcal{T}_h$ and $F$ has non-zero $d-1$-dimensional
measure.
All interior faces are collected in the set $\mathcal{F}_h^i$.
Likewise, $F$ is a boundary face if it is
the intersection of some $T^-(F)\in\mathcal{T}_h$ with $\partial\Omega$
and has non-zero $d-1$-dimensional measure.
All boundary faces make up the set $\mathcal{F}_h^{\partial\Omega} =
\mathcal{F}_h^D \cup \mathcal{F}_h^N \cup \mathcal{F}_h^O$ (subdivided into Dirichlet, Neumann
and outflow boundary faces) and we set
$\mathcal{F}_h = \mathcal{F}_h^i \cup \mathcal{F}_h^{\partial\Omega}$.
The diameter of $F\in\mathcal{F}_h$ is $h_F$ and with each $F\in\mathcal{F}_h$
we associate a unit normal vector $\vec{\nu}_F$ oriented from $T^{-}(F)$ to $T^+(F)$
in case of interior faces
and coinciding with the unit outer normal to $\Omega$ in case of boundary faces.
Every face $F$ is the image of a map $\mu_F : \hat F \to F$ with $\hat F$ the
reference element of the face.

The DG finite element space of degree $p$ on the mesh $\mathcal{T}_h$ is then
\begin{equation}
V_h^p = \left\{ v\in L^2(\Omega) : \forall\, T\in\mathcal{T}_h,
v|_T = q \circ \mu_T^{-1} \text{ with $q\in\mathbb{Q}_p^d$}\right\}
\end{equation}
where $\mathbb{Q}_p^d$ is the set of polynomials of maximum degree $p$ in dimension $d$
(all methods presented could be generalized to varying polynomial degree per direction and also per cell).
A function $v\in V_h^p$ is two-valued on an interior face $F\in\mathcal{F}_h^i$ and
by $v^-$ we denote the restriction from $T^-(F)$ and by $v^+$ the restriction
from $T^+(F)$.
For any point $\vec{x}\in F \in \mathcal{F}_h^i$ we define jump and weighted average as
\begin{align*}
\llbracket v \rrbracket (\vec{x}) &= v^-(\vec{x})-v^+(\vec{x}), &
\{ v \}_\omega (\vec{x}) &= \omega^-(\vec{x}) v^-(\vec{x}) - \omega^+(\vec{x}) v^+(\vec{x})
\end{align*}
for weights $\omega^-(\vec{x}) + \omega^+(\vec{x}) = 1$, $\omega^\pm(\vec{x}) \geq 0$.
A particular choice of the weights depending on the
diffusion coefficient $D$ has been introduced in
\cite{Ern2008,Ern01042009}:
\begin{align*}
\omega^-(\vec{x}) &= \frac{\delta_{D\nu}^+(\vec{x})}{\delta_{D\nu}^-(\vec{x}) + \delta_{D\nu}^+(\vec{x})}, &
\omega^+(\vec{x}) &= \frac{\delta_{D\nu}^-(\vec{x})}{\delta_{D\nu}^-(\vec{x}) + \delta_{D\nu}^+(\vec{x})}
\end{align*}
with $\delta_{D\nu}^{\pm}(\vec{x}) =\vec{\nu}^T(\vec{x}) D^\pm(\vec{x}) \vec{\nu}(\vec{x})$.
The definitions of jump and average are extended to boundary points
$\vec{x} \in F \in \mathcal{F}_h^{\partial\Omega}$ by
\begin{equation*}
\llbracket v \rrbracket (\vec{x}) = \{ v \}_\omega (\vec{x}) = v^-(\vec{x}) .
\end{equation*}
Finally, we denote for any domain $Q$ by
\begin{equation*}
(v,w)_Q = \int_Q v\cdot w \ dx
\end{equation*}
the $L^2$ scalar product of two (possibly vector-valued) functions, by
$|Q| = (1,1)_Q$ the measure of the set $Q$ and by
$\langle a, b \rangle = 2ab/(a+b)$ the harmonic mean of two numbers.

Following the method of lines paradigm we first discretize
equation \eqref{eq:modelproblem} in space using the WSIPG method.
This leads to the discrete in space, continuous in time problem
\begin{equation}\label{eq:discreteweakproblem}
u_h(t) \in V_h^p : \qquad \frac{d}{dt} (u_h(t),v)_\Omega + a_h(u_h(t),v;t) = l_h(v;t) \qquad \forall\, v\in V_h^p, t\in\Sigma,
\end{equation}
with the bilinear form
\begin{equation}
\label{eq:blf}
\begin{split}
a_h(u,v;t) &=
\sum_{T\in\mathcal{T}_h} \left[(D\nabla u - \vec{b}u,\nabla v)_T + (cu,v)_T\right]
+\sum_{F\in\mathcal{F}_h^{i}}
(\Phi(u^-,u^+,\vec{\nu}_F\cdot \vec{b}),\llbracket v \rrbracket)_F\\
& \quad +\sum_{F\in\mathcal{F}_h^{DO}}
(\Phi(u,0, \vec{\nu}_F\cdot \vec{b}), v)_F -\sum_{F\in\mathcal{F}_h^{iD}}
  (\vec{\nu}_F\cdot \{D\nabla u\}_\omega, \llbracket v \rrbracket)_F\\
& \quad -\sum_{F\in \mathcal{F}_h^{iD}} (\vec{\nu}_F\cdot\{D \nabla v\}_\omega , \llbracket u \rrbracket)_F
+ \sum_{F\in\mathcal{F}_h^{iD}} \gamma_{F} (\llbracket u \rrbracket , \llbracket v \rrbracket )_F
\end{split}
\end{equation}
and the right hand side functional
\begin{equation}
\begin{split}
l_h(v;t) &=  \sum_{T\in\mathcal{T}_h} (f , v)_T
- \sum_{F\in\mathcal{F}_h^{N}} (j , v)_F  \, ds
- \sum_{F\in\mathcal{F}_h^{D}} (\Phi(0,g,\vec{\nu}_F\cdot \vec{b}), v)_F  \\
&\quad -\sum_{F\in \mathcal{F}_h^{D}} (\vec{\nu}_F\cdot (D \nabla v), g)_F \, ds
+ \sum_{F\in\mathcal{F}_h^{D}} \gamma_{F} (g , v)_F.
\end{split}
\end{equation}
Here we adopted the notation $\mathcal{F}_h^{AB}=\mathcal{F}_h^{A}
\cup\mathcal{F}_h^{B}$.
The penalty factor $\gamma_{F}$ is chosen as
\begin{align*}
\gamma_{F} &= \alpha \, \langle \delta_{D\nu_F}^- , \delta_{D\nu_F}^+ \rangle M, &
M &= p (p+d-1) \frac{|F|}{\min(|T^-(F)|,|T^+(F)|)}
\end{align*}
where $\alpha$ is a free parameter typically chosen to be $\alpha=2$.
This formulation employs the upwind flux on the face which is given by
\begin{equation*}
\Phi(u^-,u^+,b_{\vec{\nu}}) = \left\{\begin{array}{ll}
u^- b_{\vec{\nu}} & b_{\vec{\nu}}\geq 0\\
u^+ b_{\vec{\nu}} & \text{else}
\end{array}\right. ,
\end{equation*}
where we denote by $b_{\vec{\nu}} = \vec{b} \cdot \vec{\nu}$ the normal flux. It is helpful to split the bilinear
form \eqref{eq:blf} into volume, interior face and boundary face contributions:
\begin{equation}
a_h(u,v;t) = \sum_{T\in\mathcal{T}_h} a_{\Omega,T}(u,v;t) + \sum_{F\in\mathcal{F}_h^{i}} a_{i,F}(u,v;t)
+ \sum_{F\in\mathcal{F}_h^{\partial\Omega}} a_{\partial\Omega,F}(u,v;t)
\end{equation}
A similar decomposition can be done for the right hand side.

Upon choosing a basis of the finite element space
$V_h^p=\text{span}\{\phi_1,\ldots,\phi_{N_h}\}$ the discrete in space problem
\eqref{eq:discreteweakproblem} is equivalent to solving
a system of ordinary differential equations
\begin{equation}
\label{eq:ODE_system}
M_h\frac{dz_h(t)}{dt} + A_h(t) z_h(t) = f_h(t)
\end{equation}
for the unknown coefficients $z_h(t)$ in the expansion
$u_h(t)=\sum_{j=1}^{N_h} (z_{h}(t))_j \phi_j$. The mass matrix $M_h$, system matrix $A_h(t)$ and the right hand
side vector $f_h(t)$ have the entries
\begin{equation*}
(M_h)_{i,j} = (\phi_j,\phi_i)_\Omega, \qquad (A_h(t))_{i,j} = a_h(\phi_j,\phi_i;t), \qquad (f_h(t))_i = l_h(\phi_i;t) \ .
\end{equation*}
A prominent advantage of DG is that the basis can be chosen such that $M_h$ is block diagonal or even diagonal
which simplifies the implementation of explicit time-stepping methods significantly.

For the time discretization of \eqref{eq:ODE_system} we use strong stability preserving explicit Runge-Kutta
methods \cite{shu:88,DiPietroErn2012}. E.g. one step of the explicit Euler method
reads
\begin{equation*}
M_h z_h^{(k+1)} = M_h z_h^{(k)} + \Delta t \left( f_h(t^{(k)}) - A_h(t^{(k)}) z_h^{(k)} \right)  .
\end{equation*}
Higher-order Runge-Kutta methods involve basically the same computations per stage.
Since we use $L_2$-orthogonal Legendre Polynomials the mass matrix $M_h$ is diagonal in our case.
Below we focus first on the matrix-free computation of $y_h = A_h z_h$ (silently dropping the time index here)
and then on the computation of a complete step of the explicit Runge-Kutta method.

\section{Sum-Factorization for Discontinuous Galerkin Methods}
\label{Sec:SpectralDG}

\subsection{Tensor Product Finite Element Functions}

On a given element $T\in\mathcal{T}_h$ a finite element function $u_h\in V_h^p$ can be expressed as
$u_h(x) = \sum_{j\in J} z_{g(T,j)} \hat\phi_j(\mu_T^{-1}(x))$ where
the $\hat\phi_j$ are the $|J|$ shape functions on the reference element $\hat T$
spanning the polynomials in $\mathbb{Q}_p^d$ and
$g: \mathcal{T}_h \times J \to \{1,\ldots,N_h\}$ is the map associating local
numbers of shape functions with numbers of global basis functions. Note that
in discontinuous Galerkin one typically choses global basis functions that have support in only one element.

The main assumption for the rest of the paper is that the shape functions have tensor
product structure,
\begin{equation}
\hat\phi_{\vec{j}}(\vec{\hat x}) = \hat\phi_{(j_1,\ldots,j_d)}(\hat x_1,\ldots,\hat x_d) =
\prod_{k=1}^d \hat\theta_{j_k}^{(k)}(\hat x_k)
\end{equation}
where $\vec{J}=J^{(1)} \times\ldots\times J^{(k)}$, $J^{(k)} = \{1,\ldots,n_k\}$, consists of
$d$-tuples enumerating the shape functions and $\hat\theta^{(k)}_{j_k}$ is
the one-dimensional basis function number $j_k$ in direction
$k$ on the reference element. In principle, the one-dimensional basis can be
different for every direction
and a different basis can be chosen for each element (anisotropic $hp$-refinement).
In our implementation, however, it is currently the same in all directions for all elements, i.e. $n_1=\ldots=n_d=n$.

For the numerical evaluation of integrals quadrature of appropriate order is used:
\begin{equation}
\int_{\hat T} f(\vec{\hat x}) \,d\hat x
= \sum_{i_1 \in I_1} \ldots \sum_{i_d \in I_d}
f\left(\xi_{i_1}^{(1)},\ldots,\xi^{(d)}_{i_d}\right) \, w_{(i_1,\ldots,i_d)} + \text{error} .
\end{equation}
Here, the quadrature formula with points $\vec{\xi_i}=\left(\xi_{i_1}^{(1)},\ldots,\xi_{i_d}^{(d)}\right)$
and weights $w_{\vec{i}} = w_{(i_1,\ldots,i_d)} = \prod_{k=1}^d w_{i_k}^{(k)}$ has tensor product form.
The $\xi_{i_k}^{(k)}$ with weights $w_{i_k}^{(k)}$ are one-dimensional quadrature points and weights
for direction $k$ with $i_k\in I^{(k)} = \{1,\ldots,m_k\}$.
$\vec{I}= I^{(1)} \times\ldots\times I^{(k)}$ is the index set of all quadrature points.
Again, quadrature order could be chosen per element and direction but in our implementation
we assume $m_1=\ldots=m_d=m$ (but \textit{not} necessarily $m=n$).

\subsection{Weak Form Evaluation}
\label{sec:weak_form_evaluation}

Evaluation of the weak form \eqref{eq:blf} requires the computation of
element-wise integrals. We now consider the part of the volume integral on element $T\in\mathcal{T}_h$ in detail for
a given test function $\phi_k$ with $k=g(T,\vec{j})$, $\vec{j}\in \vec{J}$:

\begin{equation}
\label{eq:evalblf}
\begin{aligned}
&\phantom{{}={}} (D\nabla u - \vec{b}u,\nabla\phi_k)_T\\ &= \int_T (D(\vec{x}) \nabla u - \vec{b}(\vec{x}) u) \cdot \nabla \phi_k \,dx\\
&= \int_{\hat T} \left(D(\mu_T(\vec{\hat x})) \hat S_T(\vec{\hat x}) \hat\nabla \hat u(\vec{\hat x}) - \vec{b}(\mu_T(\vec{\hat x})) \hat u(\vec{\hat
    x})\right)
\cdot \left(\hat S_T(\vec{\hat x}) \hat\nabla \hat\phi_{\vec{j}}(\vec{\hat x})\right) \Delta_T(\vec{\hat x}) \,dx\\
&\approx \sum_{\vec{i}\in I}
\left(\hat D(\vec{\xi_i}) \hat\nabla \hat u(\vec{\xi_i}) - \vec{\hat b}(\vec{\xi_i}) \hat u(\vec{\xi_i})\right) \cdot \hat\nabla\hat\phi_{\vec{j}}(\vec{\xi_i}) \Delta_T(\vec{\xi_i})  w_{\vec{i}}\\
&= \sum_{i_1\in I^{(1)}} \ldots \sum_{i_d\in I^{(d)}}
\sum_{r=1}^d \hat\partial_r\hat\phi_{\vec{j}}(\vec{\xi_i}) \left[\left(
\sum_{s=1}^d \hat D_{r,s}(\vec{\xi_i}) \hat\partial_s \hat u(\vec{\xi_i}) - \hat b_r(\vec{\xi_i}) \hat u(\vec{\xi_i})\right) \Delta_T(\vec{\xi_i})  w_{\vec{i}}\right] \\
&= \begin{aligned}[t]\sum_{r=1}^d \sum_{i_1\in I^{(1)}} \ldots \sum_{i_d\in I^{(d)}}&
\hat\partial_r \left(\hat\theta_{j_1}^{(1)}(\xi_{i_1})\cdot\ldots\cdot \hat\theta_{j_d}^{(d)}(\xi_{i_d})) \right)\\& \left[\left(
\sum_{s=1}^d \hat D_{r,s}(\vec{\xi_i}) \hat\partial_s \hat u(\vec{\xi_i}) - \hat b_r(\vec{\xi_i}) \hat u(\vec{\xi_i})\right) \Delta_T(\vec{\xi_i})
w_{\vec{i}}\right] \\ \end{aligned}\\
&= \sum_{r=1}^d \sum_{i_1\in I^{(1)}} \ldots \sum_{i_d\in I^{(d)}}
A^{(1,r)}_{j_1,i_1}\cdot\ldots\cdot A^{(d,r)}_{j_d,i_d} \ \ x^{(r)}_{(i_1,\ldots,i_d)}
\end{aligned}
\end{equation}
where now
\begin{equation}
\label{eq:gradient_matrices}
A^{(q,r)}_{\alpha,\beta} = \left\{\begin{array}{ll}
\frac{d\hat\theta_{\alpha}^{(q)}}{d\hat x}(\xi_{\beta}^{(q)})   & q=r\\[3pt]
\hat\theta_{\alpha}^{(q)}(\xi_{\beta}^{(q)}) & \text{else}
\end{array}\right., \qquad \alpha\in J^{(q)}, \beta\in I^{(q)},
\end{equation}
are small matrices which contain the values of the (derivative of) the one-dimensional basis
functions at the one-dimensional quadrature points per direction and
\begin{equation}
\label{eq:gauss_points}
x^{(r)}_{i}=\left( \sum_{s=1}^d \hat D_{r,s}(\vec{\xi_i}) \hat\partial_s \hat u(\vec{\xi_i}) - \hat b_r(\vec{\xi_i}) \hat u(\vec{\xi_i})\right) \Delta_T(\vec{\xi_i})  w_{\vec{i}}
\end{equation}
are quantities to be computed at every quadrature point involving values of (the derivatives of) the
finite element function $u$ on the element, the coefficients of the PDE and the geometry of the element $T$.

In \eqref{eq:evalblf} we use $\hat S_T(\vec{\hat x}) = (\nabla \mu_T(\vec{\hat x}))^{-T}$ in order
to transform gradients from the reference element to the real element,
set $\hat D(\vec{\hat x}) = \hat S_T^T(\vec{\hat x}) D(\mu_T(\vec{\hat x})) \hat S_T(\vec{\hat x})$,
$\hat b(\vec{\hat x}) = \hat S_T^T(\vec{\hat x}) \vec{b}(\mu_T(\vec{\hat x}))$  and abbreviate
$\Delta_T(\vec{\hat x}) = |\det\hat\nabla\mu_T(\vec{\hat x})|$.

We remark that \eqref{eq:evalblf} remains valid for nonlinear PDEs as well as high-order geometry transformations.
Only the tensor product structure of quadrature and test functions as well as linearity in the test function is used.

%
%

\subsection{Finite Element Function Evaluation}
\label{sec:fe_function_evaluation}

Equation \eqref{eq:gauss_points} requires the evaluation of a finite element function
and its gradient at all quadrature points in an element $T\in\mathcal{T}_h$.
We demonstrate the function evaluation:
\begin{equation}
\label{eq:evalfe}
\begin{split}
\hat u(\vec{\xi_i}) &= \sum_{\vec{j\in J}} z_{g(T,\vec{j})} \hat\phi_{\vec{j}}(\vec{\xi_i})\\
&= \sum_{j_1\in J^{(1)}} \ldots \sum_{j_d\in J^{(d)}} \hat\theta_{j_1}^{(1)}(\xi^{(1)}_{i_1})
    \cdot\ldots\cdot\hat\theta_{j_d}^{(d)}(\xi^{(d)}_{i_d}) \ \ z_{g(T,(j_1,\ldots,j_d))} \\
&= \sum_{j_1\in J^{(1)}} \ldots \sum_{j_d\in J^{(d)}} A^{(1)}_{i_1,j_1}\cdot\ldots\cdot A^{(d)}_{i_d,j_d} \ \ x_{(j_1,\ldots,j_d)} .
\end{split}
\end{equation}
With $A^{(q)}_{\alpha,\beta} = \theta_{\beta}^{(q)}(\xi^{(q)}_{\alpha})$, $\alpha\in I^{(q)}$, $\beta\in J^{(q)}$ and
$x_{(j_1,\ldots,j_d)}=z_{g(T,(j_1,\ldots,j_d))}$ this computation has the same structure
as that in \eqref{eq:evalblf}. The gradient can be evaluated in a similar fashion and involves
the matrices $(A^{(q,r)})^T$ from \eqref{eq:gradient_matrices}.

\subsection{Face Integral Evaluation}
\label{sec:face_integral_evaluation}

Some more notation is needed
to describe the evaluation of face integrals. The embedding of face $F$ into $T^\pm(F)$ is described by
maps $\eta_F^\pm : \hat F \to \hat T$ of the corresponding reference
elements such that $\mu_{T^\pm(F)}(\eta_F^\pm(\vec{\hat x})) = \mu_F(\vec{\hat x})$ holds.
The maps $\eta_F^\pm$ map coordinate number $q\in\{1,\ldots,d-1\}$ in $\hat F$ to coordinate number $\pi_F^\pm(q)$ in $\hat T$
and correspondingly $T^\pm(F)$.
By $q_F^\pm\in\{1,\ldots,d\}$ we denote the unique coordinate number that face $\hat F$ is perpendicular to in
the corresponding $\hat T$
and we may formally \textit{extend} the map $\pi_F^\pm$ to $\{1,\ldots,d\}$ by setting $\pi_F^\pm(d)=q_F^\pm$
(so $\pi_F^\pm$ is a permutation of $\{1,\ldots,d\}$).
Finally, since faces in the volume reference element are axi-parallel,
the component $(\eta_F^\pm(\hat x))_{\pi_F^\pm(q)} = \eta^\pm_{F,\pi_F^\pm(q)}(\hat x_q)$ is a function of one argument only
and in the case of conforming refinement it is even an isometry.

For illustration we now consider
the interior penalty term on face $F\in\mathcal{F}_h^{i}$ evaluated for a test function $\phi_k$ with support
on $T^-(F)$ and $k=g(T^-(F),\vec{j})$, $\vec{j} \in \vec{J}$.
\begin{equation}
\label{eq:eval_jump_face}
\begin{split}
(&\llbracket u \rrbracket ,\llbracket \phi_k \rrbracket)_{F} =
\int_F \llbracket u(\vec{x}) \rrbracket \phi_k(\vec{x})  \, dx
=  \int_{\hat F} \llbracket u(\mu_F(\vec{\hat x})) \rrbracket \hat\phi_{\vec{j}}(\eta_F^-(\vec{\hat x})) \Delta_F(\vec{\hat x}) \, dx \\
&\approx \sum_{i_1\in I^{(1)}} \ldots \sum_{i_{d-1}\in I^{(d-1)}}
\hat\theta_{j_{\pi_F^-(1)}}^{(\pi_F^-(1))}\left(\eta^-_{F,\pi_F^-(1)}\left(\xi^{(1)}_{i_1}\right)\right)
\cdot\ldots\cdot\hat\theta_{j_{\pi_F^-(d-1)}}^{(\pi_F^-(d-1))}\left(\eta^-_{F,\pi_F^-(d-1)}\left(\xi^{(d-1)}_{i_{d-1}}\right)\right)\\
& \qquad\qquad\hat\theta^{(q_F^-)}_{j_{q_F^-}}\left( \eta^-_{F,q_F^-}\left( 0 \right)\right)
\llbracket \hat u(\eta^-_F(\vec{\xi_i})) \rrbracket \Delta_F(\vec{\xi_i}) w_{\vec{i}} \\
&= \sum_{i_1\in I^{(1)}} \ldots \sum_{i_{d-1}\in I^{(d-1)}} \sum_{i_{d}\in I^{(d)}}
B^{(1)}_{j_{\pi_F^-(1)},i_1} \cdot\ldots\cdot B^{(d-1)}_{j_{\pi_F^-(d-1)},i_{d-1}}
B^{(d)}_{j_{q_F^-},i_d} \ x_{(i_1,\ldots,i_{d})} \ .
\end{split}
\end{equation}
where we formally introduced the one element index set $I^{(d)}=\{1\}$ and defined the matrices
\begin{equation*}
B^{(q)}_{\alpha,\beta} =
\begin{cases}
\hat\theta_{\alpha}^{(\pi_F^-(q))}\left(\eta^-_{F,\pi_F^-(q)}\left(\xi^{(q)}_{\beta}\right)\right)
& 1\leq q < d, \alpha\in J^{(\pi_F^-(q))}, \beta\in I^{(q)} \\
\hat\theta^{(q_F^-)}_{\alpha}\left( \eta^-_{F,q_F^-}\left( 0 \right)\right)
& q=d, \qquad\!\alpha\in J^{(\pi_F^-(q))}, \beta = 1
\end{cases}
\end{equation*}
as well as the coefficients
\begin{equation*}
x_i = x_{(i_1,\ldots,i_{d})} =
\llbracket \hat u(\eta^-_F(\xi^{(1)}_{i_1},\ldots,\xi^{(d-1)}_{i_{d-1}})) \rrbracket \Delta_F(\xi^{(1)}_{i_1},\ldots,\xi^{(d-1)}_{i_{d-1}}) w_{i_1,\ldots,i_{d-1}}.
\end{equation*}
Note that $\eta^-_{F,q_F^-}\left( 0 \right)\in\{0,1\}$ is independent of its argument since the face
is perpendicular to direction $q_F^-$ and therefore we can supply a dummy
argument.

Structurally we obtain the same expression as for the evaluation of the volume integral.
Similar to section \ref{sec:fe_function_evaluation} it can be shown that evaluation of the finite
element functions from both sides on a face leads to a similar expression now involving
the transposed of the matrices $B^{(q)}$.

\subsection{Sum Factorization}
\label{Sec:sumfactorization}
\eqref{eq:evalfe}
Subsections \ref{sec:weak_form_evaluation}, \ref{sec:fe_function_evaluation}
and \ref{sec:face_integral_evaluation} demonstrate that finite element function
evaluation and computation of finite element integrals 
all lead to the same abstract structure (e.g. compare to equation \eqref{eq:evalfe}):
\begin{equation}
\label{eq:Ax_component_form_2}
y_{(i_1,\ldots,i_d)}
= \sum_{j_1\in J^{(1)}}^{n_1} \ldots \sum_{j_d\in J^{(d)}}^{n_d} A_{i_1,j_1}^{(1)} \ldots A_{i_d,j_d}^{(d)}
x_{(j_1,\ldots,j_d)}
\end{equation}
for all indices $ (i_1,\ldots,i_d)\in I^{(1)}\times\ldots\times I^{(d)}$.
Assuming that $|I^{(q)}|=m$ and $|J^{(q)}|=n$ for all $q\in\{1,\ldots,d\}$ and $n/m=\rho\leq 1$
(number of quadrature points is not smaller than number of basis functions per direction)
the cost for naive evaluation of the expression \eqref{eq:Ax_component_form_2} is
\begin{equation*}
\text{Cost}_{naive}= (d+1) \rho^d m^{2d} \ \text{flops}.
\end{equation*}

A substantial reduction in the computational cost can be achieved with the
\textit{dimension by dimension} or \textit{sum factorization} algorithm
\cite{ORSZAG198070,Buis:1996:EVP:225545.225548,Melenk20014339} based
on extracting common factors in the sums:
\begin{align*}
y_{(i_1,\ldots,i_d)} &= \sum_{j_d\in J^{(d)}} \ldots \sum_{j_1\in J^{(1)}}
      A_{i_d,j_d}^{(d)} \cdot \ldots \cdot A_{i_1,j_1}^{(1)} x_{(j_1,\ldots,j_d)}\\
&= \sum_{j_d\in J^{(d)}} \ldots \sum_{j_{2}\in J^{(2)}} A_{i_d,j_d}^{(d)} \cdot \ldots \cdot A_{i_{2},j_{2}}^{(2)}
\sum_{j_d\in J^{(1)}} A_{i_1,j_1}^{(1)} x_{(j_1,\ldots,j_d)}\\
&= \sum_{j_d\in J^{(d)}} \ldots \sum_{j_{2}\in J^{(2)}} A_{i_d,j_d}^{(d)} \cdot \ldots \cdot A_{i_{2},j_{2}}^{(2)}
z^{(1)}_{(i_1,j_2,\ldots,j_{d})}\\
&= \sum_{j_d\in J^{(d)}} \ldots \sum_{j_{3}\in J^{(3)}} A_{i_d,j_d}^{(d)} \cdot \ldots \cdot A_{i_{3},j_{3}}^{(3)}
\sum_{j_{2}\in J^{(2)}} A_{i_{2},j_{2}}^{(2)} z^{(1)}_{(i_1,j_2,\ldots,j_{d})}\\
&= \sum_{j_d\in J^{(d)}} \ldots \sum_{j_{3}\in J^{(3)}} A_{i_d,j_d}^{(d)} \cdot \ldots \cdot A_{i_{3},j_{3}}^{(3)}
z^{(2)}_{(i_1,i_2,j_3,\ldots,j_{d})} = \ldots\\
&= \sum_{j_d\in J^{(d)}}  A_{i_d,j_d}^{(d)} z^{(d-1)}_{(i_{1},\ldots,i_{d-1},j_d)} = z^{(d)}_{(i_1,\ldots,i_d)} .
\end{align*}
Per direction, the \textit{sum factorization kernel}
\begin{equation}
\label{eq:sum_fact_kernel}
z^{(q)}_{(i_1,\ldots,i_q,j_{q+1},\ldots,j_{d})} = \sum_{j_q\in J^{(q)}} A_{i_q,j_q}^{(q)} z^{(q-1)}_{(i_1,\ldots,i_{q-1},j_q,\ldots,j_d)},
\qquad \forall\, i_q\in I^{(q)},
\end{equation}
needs to be carried out, where $z^{(0)} = x$. The computation in \eqref{eq:sum_fact_kernel} can be
viewed as a matrix-matrix product (albeit of rather small matrices) when $z^{(q-1)}$ is stored as a
matrix with $j_{q-1}$ as row index and the product of the other indices as column index
\cite{Buis:1996:EVP:225545.225548}. In order to facilitate recursive application, the output
$z^{(q)}$ needs to be stored such that $j_{q}$ is the new row index, which can be interpreted as a
tensor rotation of the output. If the matrix-matrix product is implemented with the columns of
$z^{q-1}$ as the slowest index and if the tensor rotation is fused into the matrix-matrix product
kernel, the resulting computation reads $z^{(q-1)}$ sequentially and writes $z^{(q)}$ sequentally,
which ensures optimal memory transfer bandwiths. Counting now the number of floating point
operations yields
\begin{equation}
\label{eq:cost_volume_sumfact}
\text{Cost}_{vol}(d,m,\rho) = 2 C_{d,\rho} m^{d+1}
\end{equation}
for volume terms and
\begin{equation}
\label{eq:cost_face_sumfact}
\text{Cost}_{face}(d,m,\rho) = 2 C_{d,\rho}  m^{d}
\end{equation}
for face terms where we set
\begin{equation*}
C_{d,\rho}  = \sum_{q=1}^{d} \rho^{q} = \left\{
\begin{array}{ll}
\rho \frac{1-\rho^d}{1-\rho} & \rho<1 \\
d & \rho=1
\end{array}\right. .
\end{equation*}
Note that in contrast to the naive evaluation all operations are fused-multiply-add operations due
to the matrix-matrix products. Achieving the optimal cost for face terms
\eqref{eq:cost_face_sumfact} requires processing the direction $q_F^\pm$ as \textit{first} direction
when evaluating a finite element function at quadrature points or as \textit{last} direction when
computing face integrals for all test functions. For moderate polynomial degrees (say $p\leq 6$) the
cost for the face terms will dominate. In that case the cost per degree of freedom will be
independent of the polynomial degree since the number of degrees of freedom per element is
$n^d=O(m^d)$.

\subsection{Matrix-free Finite Element Operator Application}
\label{sec:operator-application}

The application of the finite element operator refers to evaluate, for a given $u_h\in V_h^p$, the bilinear form $a(u_h,\phi)$ for all basis
functions $\phi$ spanning the finite element space $V_h^p$. As an extension of \cite{Kronbichler2012135} this operation can be
performed in a matrix-free setting by the following algorithm using three steps per mesh entity:
\begin{algorithm}
\caption{Matrix-free operator evaluation.}
\label{alg:AssemblyLoop}
\begin{algorithmic}
\For{$T\in\mathcal{T}_h$}
\State (1) Compute $\hat u_h, \hat\nabla\hat u_h$ for all quadrature points on corresponding $\hat T$
\State (2) Compute coefficients at all quadrature points
\State (3) Compute $a_{\Omega,T}(u_h,\phi)$ for all basis functions $\phi$ with support on $T$
\For{$F\in\mathcal{F}_h^i\cap\partial T$ and $F$ not treated yet}
\State (\phantom{1}4) Compute $\hat u_h^-, \hat u_h^+, \hat\nabla\hat u_h^-, \hat\nabla\hat u_h^+$
for all quadrature points on corresp. $\hat F$
\State (\phantom{1}5) Compute coefficients at all quadrature points
\State (\phantom{1}6) Compute $a_{i,F}(u_h,\phi^\mp)$ for all $\phi^-$, $\phi^+$ with support on $T^-(F)$, $T^+(F)$
\State (\phantom{1}7) Mark $F$ as treated
\EndFor
\For{$F\in\mathcal{F}_h^{\partial\Omega}\cap\partial T$}
\State (\phantom{1}8) Compute $\hat u_h^-, \hat\nabla\hat u_h^-$ for all quadrature points on corresponding $\hat F$
\State (\phantom{1}9) Compute coefficients at all quadrature points
\State (10) Compute $a_{\partial\Omega,F}(u_h,\phi)$ for all basis functions $\phi$ with support on $T$
\EndFor
\EndFor
\end{algorithmic}
\end{algorithm}
The computational complexity of steps (1) and (3) of the algorithm is given in \eqref{eq:cost_volume_sumfact},
the computational complexity of steps (4), (6), (8) and (10) of the algorithm is given by \eqref{eq:cost_face_sumfact}
while the computational complexity of the remaining steps (2), (5) and (9) is
proportional to the number of quadrature points.

Based on the cost of the volume and face sum factorization kernels \eqref{eq:cost_volume_sumfact}
and \eqref{eq:cost_face_sumfact} we can estimate the total cost for a full operator evaluation.
For the case $n=m=p+1$, constant full diffusion tensor, constant velocity field
and affine element transformation we obtain the following lower bound on the number of floating
point operations per degree of freedom (FLOPDOF) depending on polynomial degree $p$ and space dimension $d$:
\begin{align}
\text{FLOPDOF}_{vol}(d,p) &= 4(d+1)d(p+1) + 2 d^2 + 5d + 3, \\
\text{FLOPDOF}_{face}(d,p) &= 8(d+1)d^2 + d(2 d^2 + 8d + 18) (p+1)^{-1} .
\end{align}
Here we estimate the computational effort at each quadrature point by $2 d^2 + 5d + 3$.
The amount of data transferred (in Bytes) per degree of freedom (BDOF) can be estimated as:
\begin{align}
\text{BDOF}_{vol}(d,p) &= 16, &\text{BDOF}_{face}(d,p) &= 32 d,
\end{align}
where we assumed that all temporary variables fit into cache and only coefficients and result
need to be transfered from/to memory.
Note that face computations are more memory-intense since there are $d$ faces per element on
average and for each face we read/write data for both adjacent elements.

\begin{figure}
\begin{center}
\begin{tikzpicture}[scale=0.85]
\begin{loglogaxis}[
xlabel={Intensity in FLOPs per Byte},
ylabel={Performance in GFLOPs/s},
grid=major,
legend entries={{volume, $d=2$},{volume, $d=3$},{face, $d=2$},{face, $d=3$},roofline,matmul},
legend pos=outer north east,
]
\addplot[red,domain=1:11,samples=11,mark=*]({(4*3*2*(x+1)+2*2*2+5*2+3)/16},{min(30.4,15*(4*3*2*(x+1)+2*2*2+5*2+3)/16});
\addplot[green,domain=1:11,samples=11,mark=+]({(4*4*3*(x+1)+2*3*3+5*3+3)/16},{min(30.4,15*(4*4*3*(x+1)+2*3*3+5*3+3)/16});
\addplot[blue,domain=1:11,samples=11,mark=*]({(8*3*2*2+2*(2*2*2+8*2+18)/(x+1))/64},{min(30.4,15*(8*3*2*2+2*(2*2*2+8*2+18)/(x+1))/64});
\addplot[magenta,domain=1:11,samples=11,mark=+]({(8*4*3*3+3*(2*3*3+8*3+18)/(x+1))/96},{min(30.4,15*(8*4*3*3+3*(2*3*3+8*3+18)/(x+1))/96});
\addplot[black,domain=0.05:100,samples=100]{min(30.4,x*15)};
\addplot[black,mark=*] coordinates {(0.25,15*0.25)};
\end{loglogaxis}
\end{tikzpicture}
\end{center}
\caption{Roofline plot for one Intel Xeon E5-2698v3 2.3 GHz core with
the following hardware characteristics: $\beta=15$ GB/$s$, $\pi=30,4$ GFLOPs/$s$.}
\label{fig:roofline}
\end{figure}
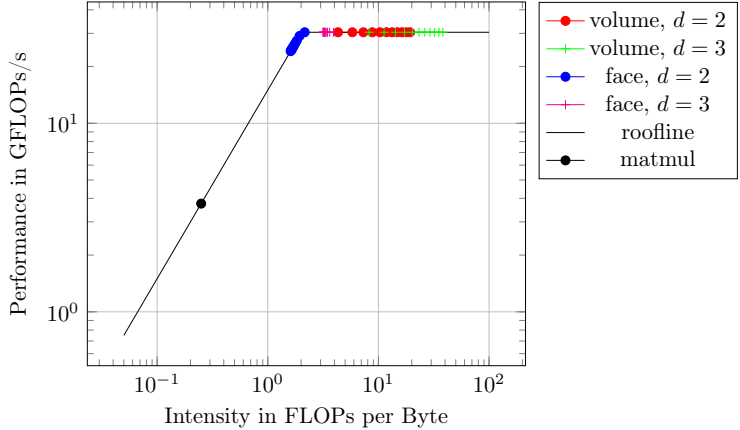

The roofline model predicts the achievable performance in GFLOPs/$s$
from the given peak performance $\pi$ (in GFLOPs/$s$) and memory bandwidth $\beta$ (in
GB/$s$) of the system as well as the theoretically determined
compute intensities
\begin{align*}
I_{vol}(d,p) &= \frac{\text{FLOPDOF}_{vol}(d,p)}{\text{BDOF}_{vol}(d,p)}, &
I_{face}(d,p) &= \frac{\text{FLOPDOF}_{face}(d,p)}{\text{BDOF}_{face}(d,p)},
\end{align*}
of the algorithm
as $$P=\min(\pi,\beta I).$$ Plotting the predicted performance $P$ over the compute
intensity $I$ gives the roofline plot in Figure \ref{fig:roofline}. It shows that the algorithm is
fully compute bound except for face terms in two space dimensions. For comparison we
include also matrix-based matrix vector multiplication which has a compute intensity of $I_{matmul}=1/4$
and performs at $3,75$ GFLOPs/$s$.

Assuming now that $3d$ computations are compute-bound for all polynomial
degrees we obtain the theoretical throughput (TPUT) in degrees of freedom (DOFs)
processed per second on a full node consisting of two Xeon E5-2698v3 2.3 GHz processors
(with 32 cores in total) as
$$\text{TPUT}(d,p) = \frac{32\pi}{\text{FLOPDOF}_{vol}(d,p)+\text{FLOPDOF}_{face}(d,p)}.$$
For $d=3$ compute intensity and throughput are shown in Figure \ref{fig:theoreticalthroughput}. Corresponding
experimental results will be shown in Section \ref{Sec:Results} below.

\begin{figure}
\begin{center}
\begin{tikzpicture}[scale=0.7]
\begin{axis}[
title={Computational Intensity},
xlabel={Polynomial Degree},
ylabel={FLOPs/DOF},
grid=major,
xmin=1,
xmax=10
]
\addplot[domain=1:10,samples=10,mark=*,]{((4*4*3*(x+1)+2*3*3+5*3+3)+(8*4*3*3+3*(2*3*3+8*3+18)/(x+1)))}; \addlegendentry{theory}
\addplot table[x=degree,y=flopsperdof-a] {benchmarks-201709/flopsperdof.txt}; \addlegendentry{measured (A)}
\end{axis}
\end{tikzpicture}
\begin{tikzpicture}[scale=0.7]
\begin{axis}[
title=Throughput,
xlabel={Polynomial Degree},
ylabel={Million DOF/sec/node},
grid=major,
xmin=1,
xmax=10
]
\addplot[domain=1:10,samples=10,mark=*,]{32*30.4*1e3/((4*4*3*(x+1)+2*3*3+5*3+3)+(8*4*3*3+2*3*(2*3*3+8*3+18)/(x+1)))}; \addlegendentry{theory}
\addplot table[x=degree,y expr=\thisrow{dofspersec-a}/1e6] {benchmarks-201709/flopsperdof.txt}; \addlegendentry{measured (A)}
\end{axis}
\end{tikzpicture}
\end{center}
\caption{Theoretical compute intensity and throughput for a $3d$ computation
on an axiparallel mesh with constant diffusion tensor assuming $\pi=32\cdot 30,4$ GFLOPs/s for the full node.}
\label{fig:theoreticalthroughput}
\end{figure}
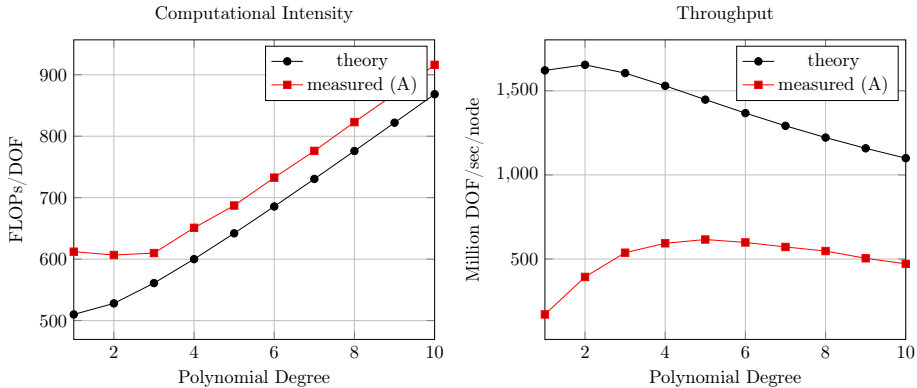

\section{Implementation Issues}
\label{Sec:Implementation}

We have realized the numerical scheme outlined in the previous section within
the PDELab \cite{pdelabalgoritmy} finite element framework which is based
on the Dune software framework \cite{Dune2008b}. Within the
EXA-DUNE project \cite{SppExaSymposiumHardware} the Dune software framework is
currently being prepared for current and future HPC architectures.
By design, Dune and PDELab are very general and allow for different mesh types and general
discretizations. This flexibility sacrifices, to some extent,
performance for generality and in order to obtain a high-performance implementation
of sum factorized DG assembly, we had to redesign parts of the underlying framework.
In this Section, we will focus on two different aspects that were crucial for the performance of our code: the vectorization
strategy to exploit SIMD parallelism of modern CPU architectures and the block structure of DG vectors that can be used to
reduce the amount of memory transfers.

Peak performance of modern many-core CPUs can only be achieved with SIMD vector instructions.
Here we are concerned with Intel Haswell processors offering the AVX2 instruction set which operates
on four double-precision floating point numbers in parallel. Although matrix-matrix multiplication
of sufficiently large matrices vectorizes very well, the small size of the matrices involved in
sum factorization is a challenge. For example, a single step of the sum-factorization algorithm
given by equation \eqref{eq:sum_fact_kernel} in
three dimensions for polynomial degree three involves the product of a $4\times 4$ and
a $4\times 16$ matrix. While the performance for multiplication of small matrices is
meanwhile addressed on a library level by \texttt{libxsmm} \cite{libxsmm_2016} we provide our own
implementation based on the following principles:
\begin{enumerate}[i)]
\item Increase size of the matrices by collecting several sum factorization steps. The FMA instructions in
current Intel processors have a latency of 4-5 processor cycles and
the processor can issue two such instructions per cycle. Therefore one requires 8-10
independent summation chains of 4 values each to fully exploit the floating
point performance of the processor.
\item Aim at at least one matrix dimension being a multiple of the SIMD vector unit (four here).
\item Arrange the output matrix of the matrix-matrix product in a way that enables recursive application of the
steps \eqref{eq:sum_fact_kernel} and avoid a seperate tensor rotation. This results in a non-standard
matrix-matrix product which is currently not available in the \texttt{libxsmm} library at the time of writing.
\end{enumerate}
There are several ways to increase the amount of work as stated in principle i).
One option is to work simultaneously on $k$ cells or facets where $k$ is the SIMD width ($k=4$ in our case).
This has the advantage of producing an innermost loop of size $k$ and vectorization is independent of
polynomial degree, quadrature order and spatial dimension. On the other hand the data of $k$ cells needs to be accessed
simultaneously which puts more pressure on the cache for larger polynomial degrees.

Another option, which we have chosen in this work, is to exploit the fact that the finite element method applied to second order PDEs
involves functions as well as their gradients and to work on just one element or facet at a time.
This is particularly appealing in three spatial dimensions for AVX2, where we can
group together the four values $[\partial_{x_1} u(\vec{\xi}),\partial_{x_2} u(\vec{\xi}),\partial_{x_3} u(\vec{\xi}), u(\vec{\xi})]$ at one quadrature point $\vec{\xi}$
into a single SIMD register.

Besides sum factorization, matrix-free operator evaluation
involves computations at each quadrature point in steps (2), (5) and (9) of algorithm \ref{alg:AssemblyLoop} which
need to be vectorized as well. As the computations at different quadrature points
are independent we process $k$ quadrature points in parallel using SIMD instructions.
This requires an additional in-place transpose operation illustrated here for $k=4$:
\begin{gather*}
\begin{array}{cllllc}
 [ & \partial_{x_1} u(\vec{\xi}_l), & \partial_{x_2} u(\vec{\xi}_l), & \partial_{x_3} u(\vec{\xi}_l), & u(\vec{\xi}_l) & ] \\
{[} & \partial_{x_1} u(\vec{\xi}_{l+1}), & \partial_{x_2} u(\vec{\xi}_{l+1}), & \partial_{x_3} u(\vec{\xi}_{l+1}), & u(\vec{\xi}_{l+1}) & ] \\
{[} & \partial_{x_1} u(\vec{\xi}_{l+2}), & \partial_{x_2} u(\vec{\xi}_{l+2}), & \partial_{x_3} u(\vec{\xi}_{l+2}), & u(\vec{\xi}_{l+2}) & ] \\
{[} & \partial_{x_1} u(\vec{\xi}_{l+3}), & \partial_{x_2} u(\vec{\xi}_{l+3}), & \partial_{x_3} u(\vec{\xi}_{l+3}), & u(\vec{\xi}_{l+3}) & ] \\
\end{array}\\[.5em]
\downarrow\\[.5em]
\begin{array}{cllllc}
{[} & \partial_{x_1} u(\vec{\xi}_l), & \partial_{x_1} u(\vec{\xi}_{l+1}), & \partial_{x_1} u(\vec{\xi}_{l+2}), & \partial_{x_1} u(\vec{\xi}_{l+3}) & ] \\
{[} & \partial_{x_2} u(\vec{\xi}_l), & \partial_{x_2} u(\vec{\xi}_{l+1}), & \partial_{x_2} u(\vec{\xi}_{l+2}), & \partial_{x_2} u(\vec{\xi}_{l+3}) & ] \\
{[} & \partial_{x_3} u(\vec{\xi}_l), & \partial_{x_3} u(\vec{\xi}_{l+1}), & \partial_{x_3} u(\vec{\xi}_{l+2}), & \partial_{x_3} u(\vec{\xi}_{l+3}) & ] \\
{[} & u(\vec{\xi}_l), & u(\vec{\xi}_{l+1}), & u(\vec{\xi}_{l+2}), & u(\vec{\xi}_{l+3}) & ] \\
\end{array}
\end{gather*}
A second transpose operation is necessary after the quadrature point computations
to achieve the correct layout for the subsequent
sum factorization step.
While the transpose steps add an additional overhead, this overhead is very small and can mostly be hidden behind the actual
computations by loop unrolling. The remaining overhead is outweighed by the performance advantage in the sum factorization steps.

Our implementation precomputes and caches the local coordinates and associated weights of all quadrature points. These are then streamed
from memory during steps (2), (5) and (9) of Algorithm~\ref{alg:AssemblyLoop}, which allows us to use a flat iteration space during those
steps. In our experience, the advantage of not having to use nested iteration or reconstruct the multi index $\vec{i}$ from a flat index by
far outweighs the increased L1 cache pressure of having to load an additional $d+1$ values per quadrature point.

We have implemented separate code paths for axis-parallel, affine and multilinear geometries.
In the case of multilinear geometries the transformation $x=\mu_T(\hat{x})$ from the reference
element $\hat{T}$ to the real element $T$ is a $d$-valued $Q_1$ finite element function
and sum-factorization can be used to evaluate its values and derivatives at quadrature points.
The same approach could also be used for higher order geometries.

The evaluation of face integrals poses an additional challenge:
As outlined in Section \ref{Sec:sumfactorization}, it is important to start (respectively end)
with the face normal direction $q_F^\pm$ in steps (4), (8) (respectively steps (6) and (10)) of algorithm \ref{alg:AssemblyLoop}
to obtain optimal computational complexity for sum factorization in face integrals.
This does, however, complicate the memory access patterns in the first (respectively last) step of the
sum factorization kernel, which causes a substantial performance reduction
which is further exacerbated by the reduced size of the involved matrix (compared to volume integrals).
In order to minimize this performance impact, we use C++ metaprogramming to
generate dedicated versions of the face sum factorization kernels for all combinations of
face normal directions.

Like other general-purpose finite element discretization frameworks,
PDELab by default uses a temporary buffer to gather all degrees of freedom associated with
a single element contiguously into memory before processing the element. Likewise, the computational results associated
with one element are stored in a buffer and then scattered to corresponding locations in the
global data structure. These gather and scatter operations are, however, not necessary
in discontinuous Galerkin methods where the global data may already be arranged contiguously in memory.
We thus extended PDELab with a DG-specific code path that exploits the block structure of DG problems
avoiding any superfluous copy operations.
%

%
%
%

\section{Numerical Results}
\label{Sec:Results}


\begin{table}
  \centering
  \footnotesize
\begin{tabularx}{\textwidth}{cYGYGYFYG}
\toprule
& \multicolumn{2}{c}{Matrix-free A}& \multicolumn{2}{c}{Matrix-free B} & \multicolumn{2}{c}{Matrix-based} & \multicolumn{2}{c}{Matrix Assembly} \\ \cmidrule(lr){2-3} \cmidrule(lr){4-5} \cmidrule(lr){6-7} \cmidrule(lr){8-9}
$p$ & \multicolumn{1}{c}{\si[per-mode=fraction]{\dof\per\second}} & \multicolumn{1}{c}{\si[per-mode=fraction]{\giga\flop\per\second}} & \multicolumn{1}{c}{\si[per-mode=fraction]{\dof\per\second}} & \multicolumn{1}{c}{\si[per-mode=fraction]{\giga\flop\per\second}} & \multicolumn{1}{c}{\si[per-mode=fraction]{\dof\per\second}} & \multicolumn{1}{c}{\si[per-mode=fraction]{\giga\flop\per\second}} & \multicolumn{1}{c}{\si[per-mode=fraction]{\dof\per\second}} & \multicolumn{1}{c}{\si[per-mode=fraction]{\giga\flop\per\second}} \\
\specialrule{\lightrulewidth}{\belowrulesep}{\belowrulesep}
1 & 170158953.810 & 104.158182879 & 118979587.729 & 321.037185434 & 206451612.9 & 24.3316669935 & 16039016.302 & 344.537537152 \\
2 & 392974786.299 & 238.415477889 & 252483269.110 & 449.598013872 & 64200044.85 & 27.2839491589 & 8709205.2273 & 371.308869533 \\
3 & 537750267.971 & 327.958675764 & 330453991.322 & 523.700278739 & 26851310.34 & 29.8372413793 & 4662890.30518 & 368.195374202 \\
4 & 595014243.273 & 387.318881958 & 388214296.397 & 559.995679261 & 9544347.82 & 23.4740484429 & 2574016.32326 & 300.728647745 \\
5 & 617156313.613 & 424.131804157 & 403296404.441 & 568.113051526 & 4584939.03 & 23.3417548171 & 1925972.94389 & 307.258859483 \\
6 & 598962766.081 & 438.799707698 & 406075144.913 & 562.979680371 & 2310736.84 & 23.4839329435 & 1210878.64306 & 230.908894365 \\
7 & 570054393.153 & 442.435533605 & 398006289.990 & 555.586463942 & 1462449.09 & 30.3253444017 & 664872.60321 & 142.838982665 \\
8 & 540645570.859 & 444.910841687 & 384895618.842 & 541.184528147 & 614493.65 & 24.0574267317 & 873503.72048 & 197.904571983 \\
9 & 504835956.761 & 438.602901953 & 370258969.165 & 529.864853667 & 297674.41 & 26.1953488372 & 700825.043546 & 132.758132756 \\
10 & 471259916.574 & 431.763201629 & 359375741.589 & 523.586101106 & \multicolumn{1}{c}{---} & \multicolumn{1}{c}{---} & \multicolumn{1}{c}{---} & \multicolumn{1}{c}{---} \\
\bottomrule
\end{tabularx}
\caption{Full operator application, 2x Intel Xeon E5-2698v3 2.3 GHz, all 32 cores busy for two
  problems of different complexity and for comparison matrix assembly for the simpler problem and
  matrix-based operator application. Note that the matrix-based computations use significantly
  smaller problem sizes due to memory constraints.}
\label{fig:throughput-table}
\end{table}

In order to measure the efficiency of our numerical approach and implementation, we are most importantly interested in the time to
solution for a given accuracy. As the aim of this paper is high-performance operator application, we are right now merely interested
in the time required for this operation. In order to keep operator applications comparable across different polynomial degrees with
different problem sizes and matrix-free / matrix-based approaches, we measure the number of DOFs/s that we can process during the
operator application.

Moreover, as the operator application itself is a completely local operation and does not require communication, we begin by studying its
behavior on a single node compute node. For all of the following measurements, we used a node of our workgroup cluster that is equipped with
2x Intel Xeon E5-2698v3 2.3 GHz Haswell processors (16 cores and 40 MiB of L3 cache each) and 128 GiB of DDR3-2133 RAM with deactivated
Hyperthreading and deactivated Turbo Mode. The peak floating point performance of this processor can only be attained using fused
multiply-add instructions; at the maximum throughput rate of 2 FMA instructions/cycle, this processor is capable of a theoretical peak of
486.4 GFLOPs (while the vector math units are active, the processor speed is reduced to 1.9 GHz for thermal management), yielding a total of
972.8 GFLOPs / node. In order to achieve this performance, the calculations must expose a considerable amount of independent operations, as
the FMA instructions on this processor have a latency of 5 clock cycles and we thus need to expose $2*5*4 = 40$ independent FMA chains.
Moreover, as the processor relies on SIMD registers of width 4, the operations must be grouped into blocks of that size operating on
adjacent memory locations as this processor does not support strided loads and stores with good performance. This problem is further
exacerbated by the limited number of available processor registers (16), which requires us to reuse registers when feeding data to the FMA
units of the processor, as each of the 10 vector operations in flight at peak performance requires 3 registers (two input registers and one
input / output register).

In order to avoid spurious performance effects due to under-utilized resources, we run our benchmarks as MPI-parallel programs with 32
ranks, one pinned to each core of the node. The meshes for problems of different polynomial degree are sized in such a way that the number
of DOFs is approximately constant ($\approx 4\cdot 10^8$), resulting in input and output vectors of about 3.2 GiB each in total to
avoid any possible cache effects.

\subsection{Model Problems}

As we want to model the performance as part of a much larger computation involving additional nodes, we set up a test problem with periodic
boundary conditions in all directions. PDELab implements periodic boundaries in the same way as regular processor boundaries, which results
in an ideally load-balanced computation.

We test our implementation with three different problems of increasing computational intensity. All problems are
convection-diffusion-reaction problems with tensor-valued diffusion coefficients and use weighted averaging across cell boundaries for
discontinuous coefficients:
\begin{description}
\item[Problem A] has axis-parallel cells and coefficients that are constant per cell. It is a baseline model that we expect to perform best
  in terms of throughput (DOFs/sec), but not necessarily in efficiency (GFLOPs/sec), due to the low amount of work per DOF.
\item[Problem B] is an intermediate problem with affine geometries and coefficients that are polynomials in the coordinates.
\item[Problem C] is a computationally intensive problem with multilinear geometries and coefficients that are polynomials in the
coordinates.
\end{description}

A central question related to these problems is the required integration order $q$: For Problem A, we can choose $q = 2p$ (which yields
$n = m = p+1$ as in the theoretical discussion in \ref{sec:operator-application}). As the diffusion tensor in Problem B/C is a polynomial in
$\vec{x}$, we increase $q$ to $q = 2p + 4$. Problem C models a situation where most/all geometries are nonlinear and as such, the inverse
determinant of the geometry Jacobian becomes a rational function in $\vec{x}$. In contrast to first-order PDEs, where this term cancels out, it is
very much present in our second-order example. In order to minimize the quadrature error, we have (somewhat arbitrarily) picked $q = 3p + 4$
for the final problem. Note that in actual applications, non-linear geometries can often be restricted to small parts of the mesh. By
constrast, we have chosen Problem C to highlight the performance impact of more complex geometries and/or coefficients.

\subsection{Comparison To Matrix-Based Operator Application}

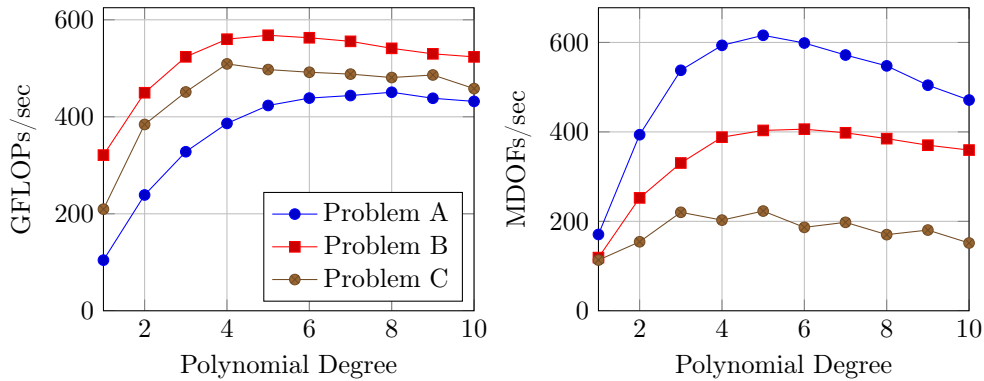
\begin{figure}[tbh]
  \centering
  \begin{tikzpicture}
    \begin{axis}[
      name=plot1,
      width=.5\textwidth,
      xlabel={Polynomial Degree},
      ylabel={GFLOPs/sec},
      ymin=0,
      xmin=1,
      xmax=10,
      grid=major,
      legend entries={{Problem A},{Problem B},{Problem C}},
      legend pos=south east
      ]
      \addplot table[x=degree,y=flopspersec-a] {benchmarks-201709/performance.txt};
      \addplot table[x=degree,y=flopspersec-b] {benchmarks-201709/performance.txt};
      \addplot table[x=degree,y=flopspersec-c] {benchmarks-201709/performance.txt};
    \end{axis}
    \begin{axis}[
      at=(plot1.right of south east), anchor=left of south west,
      width=.5\textwidth,
      xlabel={Polynomial Degree},
      ylabel={MDOFs/sec},
      ymin=0,
      xmin=1,
      xmax=10,
      grid=major,
      ]
      \addplot table[x=degree,y expr=\thisrow{dofspersec-a}/1e6] {benchmarks-201709/performance.txt};
      \addplot table[x=degree,y expr=\thisrow{dofspersec-b}/1e6] {benchmarks-201709/performance.txt};
      \addplot table[x=degree,y expr=\thisrow{dofspersec-c}/1e6] {benchmarks-201709/performance.txt};
    \end{axis}
  \end{tikzpicture}
  \caption{Floating point performance in GFLOPs/sec and throughput in MDOFs/sec for full operator application, 2x Intel Xeon E5-2698v3 2.3
    GHz for all model problems}
  \label{fig:problem-comparison}
\end{figure}

Table~\ref{fig:throughput-table} shows the throughput and the hardware efficiency of our matrix-free code for Problem A and Problem B. For
comparison, we also assemble the matrix of the operator (again using our sum-factorized implementation) and apply the resulting matrix to an
input vector using matrix-vector multiplication. Those results are also shown in Table \ref{fig:throughput-table}; they were obtained using
PETSc \cite{petsc-user-ref,petsc-web-page} using block compressed row format to exploit the block structure of DG. Note that the problems
had to be scaled down a lot for the matrix-based computations in order to fit the matrix into memory: For $p=3$, the matrix-based problem is
about 15 times smaller, while for $p=9$, the factor is $\approx 1300$. Our matrix-free implementation outperforms the matrix-vector product for
$p>1$ or $p>2$ (simple/expensive problem). However, even at low degrees the performance advantage of the matrix-based version is somewhat
reduced by the fact that assembling the matrix also takes a significant amount of time in addition to severely reducing the number of DOFs
that can fit into a given amount of RAM.

Figure~\ref{fig:problem-comparison} compares throughput and floating point performance of our implementation for the three problems
described above. With increasing complexity (i.e., more work per DOF), throughput decreases by a factor of $\approx3$ when going from Problem A to
Problem C. On the other hand, the higher computational intensity affords better  as demonstrated by the the GFLOPs/sec rates,
which reach up to about 60\% peak for Problem B. Problem C's performance falls between that of the other two problems; while it provides the
most work/DOF, the additional sum factorization operations for the coordinates are only polynomials with $p_C=1$, and with increasing $p$ we
start to exhaust the L1 cache due to the large number of quadrature points.


\subsection{Intra-Node Scalability}
\label{sec:intra-node-scalability}

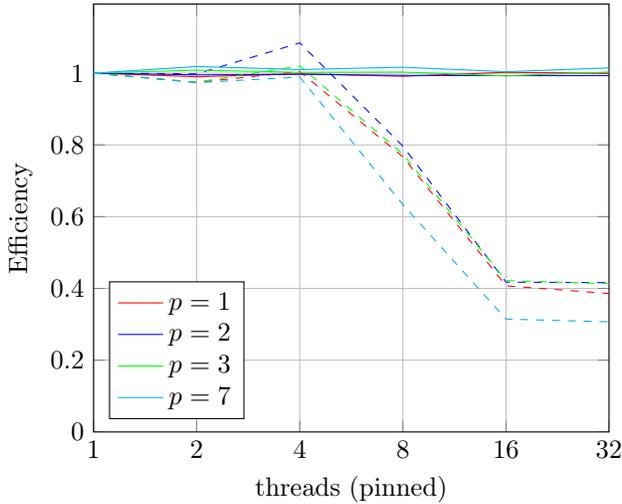
\begin{figure}[tbh]
  \centering
  \begin{tikzpicture}
  \tikzset{every mark/.append style={scale=.5}}
    \begin{semilogxaxis}[
      legend columns=1,
      legend cell align=left,
      legend entries={$p=1$,$p=2$,$p=3$,{$p=7$}},
      legend pos=south west,
      ymin=0,
      xlabel=threads (pinned),
      ylabel=Efficiency,
      xmin=1,xmax=32,
      xtick={1,2,4,8,16,32},
      xticklabels={1,2,4,8,16,32},
      grid=major,
      ]

      \addplot[red] table[x index=0,y index =1] {scaling2.txt};
      \addplot[red,dashed,forget plot] table[x index=0,y index =1] {mat-scaling.txt};
      \addplot[blue] table[x index=0,y index =3] {scaling2.txt};
      \addplot[blue,dashed,forget plot] table[x index=0,y index =2] {mat-scaling.txt};
      \addplot[green] table[x index=0,y index =5] {scaling2.txt};
      \addplot[green,dashed,forget plot] table[x index=0,y index =3] {mat-scaling.txt};
      \addplot[cyan] table[x index=0,y index =7] {scaling2.txt};
      \addplot[cyan,dashed,forget plot] table[x index=0,y index =4] {mat-scaling.txt};

    \end{semilogxaxis}
  \end{tikzpicture}\\[-1em]
  \caption{Weak scaling of operator application on single node. Processes pinned to individual cores, runs with 1-16 processes
    pinned to a single socket and 32 processes across two sockets. Solid lines denote matrix-free computations, dashed lines matrix-based.
  }
  \label{fig:intranode-scaling}
\end{figure}

Next, we look at the weak scalability of matrix-free and matrix-based operator application on a single compute node, which provides insight
into the load placed on parts of the processor shared across multiple cores (memory interfaces, I/O, etc.). As shown in
Figure~\ref{fig:intranode-scaling}, our matrix-free operator application scales perfectly to all cores of the machine, while the
matrix-based version breaks down beyond 4 cores per socket. This behavior can be explained by the roofline model
(cf.~Figure~\ref{fig:roofline}): The matrix-based implementation is memory-bound, and as the processor in our machine is equipped with 4
memory controllers, their bandwidth has to be shared across multiple cores. The matrix-free code is compute-bound and thus not restricted by
resources shared between multiple cores. The cores are pinned to sockets in such a way that we first saturate one socket before allocating
cores on the second socket. Thus the step from 16 to 32 cores is balanced in the matrix-based case by the additional 4 memory controllers on
the second socket.

\subsection{Realistic Example Problem}
\label{sec:example-problem}

\begin{figure}
  \centering
  \includegraphics[width=.8\textwidth]{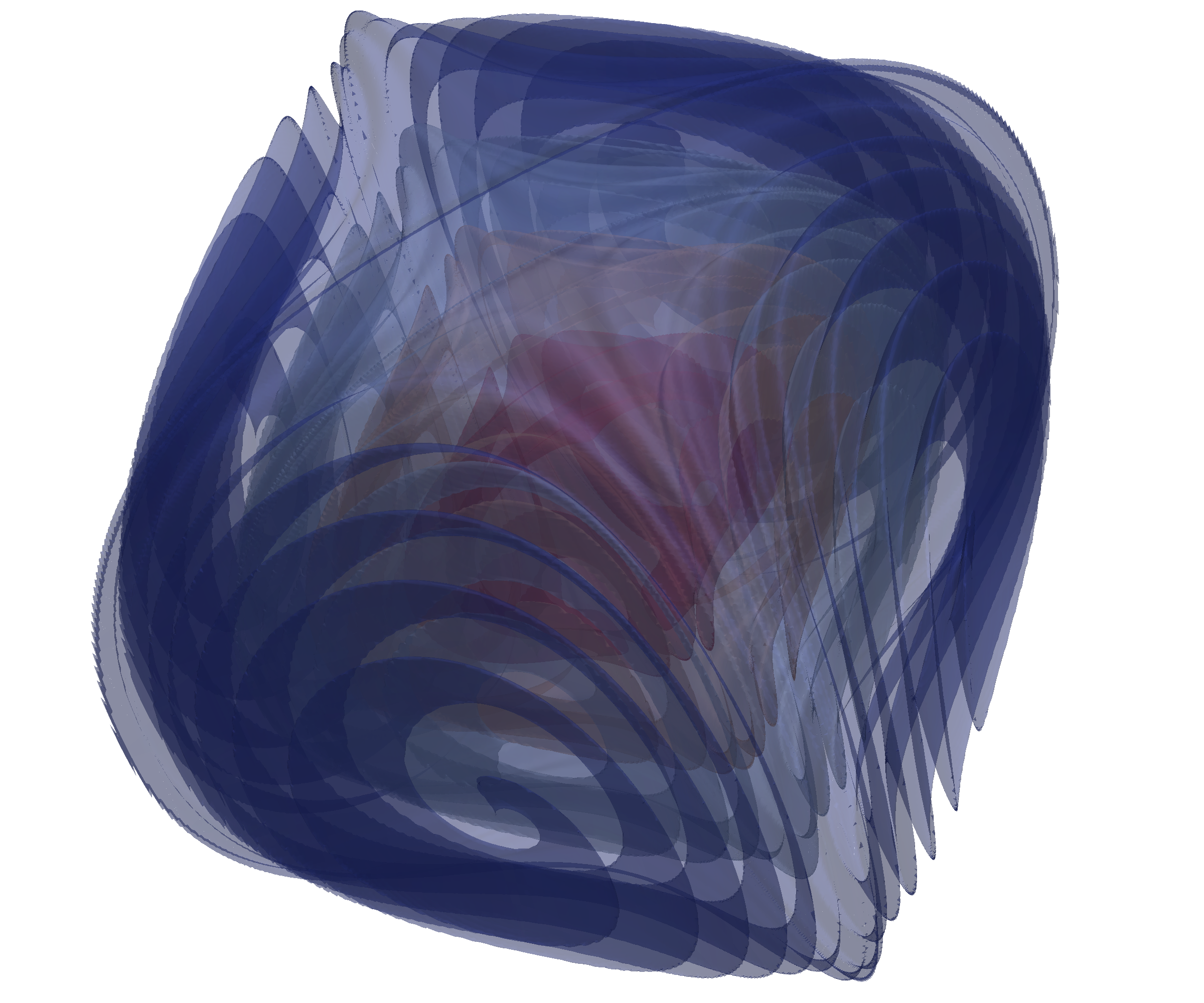}
  \caption{Isocontours for $u \in \{0.13, 0.25, 0.38, 0.51, 0.64\}$ at $t=140$ for the Taylor-Green example problem with
    polynomial degree 5 and mesh size $240^3$ ($\approx 3\cdot10^9$ DOFs, 165000 time steps, 143000 core hours
  on 6144 cores including I/O).}
  \label{fig:taylorgreen}
\end{figure}

In addition to measuring the performance of our code in isolation, we also apply it to a more complete problem setting. We study the full
instationary transport problem \eqref{eq:modelproblem} with highly convection-dominated flow on a periodic domain $\Omega = (-\pi,\pi)^3$ with
$D = 5\cdot10^{-6}\cdot\mathbb{I}$, $c = 0$, $f = 0$ and a fixed flow field of Taylor-Green type \cite{Taylor1937}:
\begin{equation}
  \label{eq:2}
  \vec{b}(\vec{x}) =
    \begin{pmatrix}
      \cos(x_1) \sin(-x_2) \sin(-x_3)\\[.2em]
      \frac{1}{2} \sin(x_1) \cos(-x_2) \sin(-x_3)\\[.2em]
      \frac{1}{2} \sin(x_1) \sin(-x_2) \cos(-x_3)\\
    \end{pmatrix}.
\end{equation}
We then compute the evolution of an initial concentration $u_0$ given by a Gaussian centered around $x_0 = (\frac{1}{2},0,0)$ as
$u_0 = \exp(-2\|x-x_0\|^2)$.

This problem is computationally quite expensive due to the trigonometric function evaluations at each quadrature point and the increased
quadrature order required for those functions (we have chosen $q = 2*p + 4$ here). The trigonometric functions are evaluated using optimized
and vectorized implementations from \cite{Fog:Vectorclass}.

In order to limit the scope of this discussion, we use a strong stability preserving second-order explicit time stepping scheme (Heun). In
combination with the DG discretization, this yields a block-diagonal mass matrix, and by choosing Legendre basis functions, we even obtain a
point-diagonal mass matrix, eliminating any occurrences of non-trivial linear equation systems in our solver. We calculate the inverse mass
matrix in a setup step, and as a result, our solver only performs matrix-free residual calculations and products of a diagonal matrix
with a vector, combined with halo exchanges to keep the solution consistent across subdomains.

The periodic boundary conditions of the problem ensure that all subdomains of the problem are identically structured, which guarantees ideal
load balancing.

Figure~\ref{fig:taylorgreen} shows the concentration at $t=140$ for a simulation with $240^3$ cells and polynomial degree
$p=5$. As can be seen in the image, the prescribed flow field creates a very complex structure that requires a very high resolution to
capture the solution well.

\subsection{Scalability Studies for Realistic Example Problem}
\label{sec:scalability-studies}

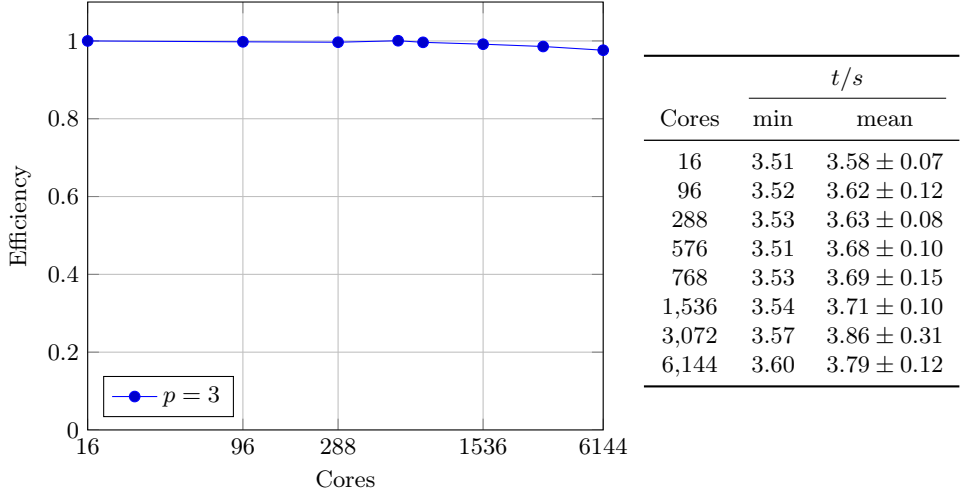
\begin{figure}[tbh]
  \centering
  \small
  \begin{tikzpicture}
    [baseline={(0,2.7)}]
    \begin{semilogxaxis}[
      ymin=0,
      ymax=1.1,
      grid=major,
      xmin=16,
      xmax=6144,
      xtick={16,96,288,1536,6144},
      xticklabels={16,96,288,1536,6144},
      legend entries={{$p=3$}},
      legend pos=south west,
      xlabel={Cores},
      ylabel={Efficiency},
      ]
      \addplot table[x expr={\thisrow{node_count}*16},y expr={3.514794921875/\thisrow{bench/halo-blocked}}] {benchmarks-201709/weak-scalability.txt};
    \end{semilogxaxis}
  \end{tikzpicture}
  \begin{filecontents*}{benchmarks-201709/weak-scalability-table.txt}
cores	min	mean
16.0	3.51479492188	3.58 +- 0.07
96.0	3.52280555556	3.62 +- 0.12
288.0	3.52641520894	3.63 +- 0.08
576.0	3.51311156836	3.68 +- 0.10
768.0	3.52746111111	3.69 +- 0.15
1536.0	3.54461162551	3.71 +- 0.10
3072.0	3.56589084201	3.86 +- 0.31
6144.0	3.60080555556	3.79 +- 0.12
\end{filecontents*}
\pgfplotstabletypeset[
    col sep=tab,
    columns={cores,min,mean},
    columns/min/.style={
      column name=$\text{min}$,
      precision=2,
      zerofill,
    },
    columns/mean/.style={
      column name=$\text{mean}$,
      string type,
      column type={S[
        table-format=1.2,
        round-mode=figures,
        separate-uncertainty,
        table-align-uncertainty,
        table-figures-uncertainty=3,
        table-auto-round
        ]},
    },
    columns/cores/.style={
      column name=Cores,
      int detect,
    },
    every head row/.style={
      before row={
        \toprule
        & \multicolumn{2}{c}{$t/s$}\\
        \cmidrule(lr){2-3}
      },
      after row=\midrule,
    },
    every last row/.style={
      after row=\bottomrule
    },
    ]
    {benchmarks-201709/weak-scalability-table.txt}
    \caption{Efficiency and run times for weak scalability on IWR compute cluster (416 nodes with 2 x E5-2630 v3 each, 64 GiB / node, QDR
      Infiniband). The plot is based on the fastest times, mean values show a large amount of timing jitter as seen in the table.}
  \label{fig:weak-scaling}
\end{figure}

\begin{figure}[tbh]
  \centering
  \begin{tikzpicture}
    \begin{loglogaxis}[
      grid=major,
      xmin=1,
      xmax=6144,
      ylabel shift=-.4em,
      xtick={1,12,96,768,6144},
      xticklabels={1,12,96,768,6144},
      legend entries={{$p=1$},{$p=3$},{$p=5$},{$p=7$}},
      legend pos=north east,
      xlabel={Cores},
      ylabel={$t/s$},
      ]
      \addplot table[x=degree,y=Q1] {benchmarks-201709/strong-scalability.txt};
      \addplot table[x=degree,y=Q3] {benchmarks-201709/strong-scalability.txt};
      \addplot table[x=degree,y=Q5] {benchmarks-201709/strong-scalability.txt};
      \addplot table[x=degree,y=Q7] {benchmarks-201709/strong-scalability.txt};
      \addplot [domain=1:6144] {80/x};
    \end{loglogaxis}
  \end{tikzpicture}
  \caption{Run times for strong scalability on IWR compute cluster (416 nodes with 2 x E5-2630 v3 each, 64 GiB / node, QDR Infiniband)}
  \label{fig:strong-scaling-fig}
\end{figure}
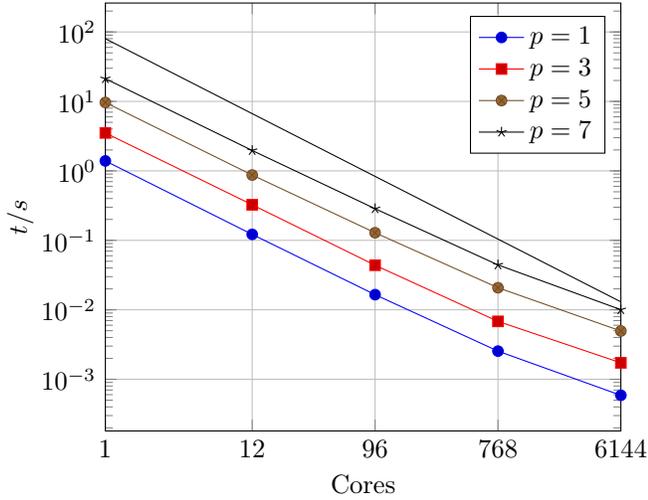

As we are ultimately interested in HPC, we investigate the scalability of our implementation on the moderately sized \emph{bwfordev}
development cluster in Heidelberg, which consists of 416 compute nodes with 2 x E5-2630 v3 processors each (Haswell architecture with 16
cores~/~node, 64 GiB~/~node) and uses a QDR Infiniband interconnect with a fully connected two-level topology. Figure~\ref{fig:weak-scaling}
shows the efficiency of weakly scaling our realistic example problem \ref{sec:example-problem} for $p=3$ with $\approx 10^6$ DOFs~/~core from 1 to 384 nodes of the cluster.

Our solver is only made up of local computations and halo exchanges, so we expect the mostly flat efficiency curve observed in our
measurements; the slight performance degradation at larger problem sizes is mostly due to the increased communication jitter as the program
starts to span more of the cluster's communication infrastructure.

The results of our strong scalability benchmarks are shown in Figure~\ref{fig:strong-scaling-fig}. Due to the low memory footprint of our
matrix-free solver, we are able to measure scalability from 1 to 6144 cores on up to 384 nodes. Across this range, the mesh size per core
shrinks from $48^3$ to $3\times3\times2$ cells, which corresponds to between $10^2$ and $10^5$ DOFs per core at 6144 cores, depending on $p$. For these
very small working sets, scalability suffers mostly because our implementation is currently not able to overlap computation and halo
communication, which can be seen in more detail in Table~\ref{fig:strong-scaling-tab}.

\begin{table}[h]
  \centering
  \small
  \begin{filecontents*}{benchmarks-201709/strong-scalability-table.txt}
cores	time_1	eff_1	time_2	eff_2	time_3	eff_3	time_4	eff_4	time_5	eff_5	time_6	eff_6	time_7	eff_7	time_8	eff_8
1.0	1.39091	1.0	2.63411	1.0	3.51836	1.0	6.16352	1.0	9.65249	1.0	14.4547	1.0	21.2497	1.0	31.3775	1.0
12.0	0.121871	0.951080787609	0.237456	0.924420383847	0.325015	0.902101954269	0.54912	0.935363247863	0.87113	0.92336868971	1.30079	0.926020597739	1.96718	0.900176055741	2.888	0.905398776547
96.0	0.0165197	0.877052599825	0.0316847	0.865990393891	0.0438419	0.835948791757	0.0782711	0.820268698579	0.12861	0.781795901045	0.194538	0.773986530481	0.285402	0.775576350785	0.42067	0.776972349664
768.0	0.00254124	0.712675988559	0.00509011	0.673822516442	0.00681274	0.672445729129	0.0120792	0.664399684306	0.0207832	0.604735861377	0.0311888	0.603460984659	0.0442931	0.624676985994	0.0657839	0.621065637514
6144.0	0.000587705	0.385201914474	0.00121699	0.352286248158	0.00172599	0.331780450399	0.00253955	0.395021591752	0.0049536	0.317151827816	0.00647358	0.363423792522	0.0100378	0.344558571205	0.0147238	0.346854410815
\end{filecontents*}
\pgfplotstabletypeset[
    columns={cores,time_1,eff_1,time_3,eff_3,time_5,eff_5,time_7,eff_7},
    every even column/.style={
      column name=$\eta$,
      precision=2,
      zerofill,
    },
    every odd column/.style={
      column name=$t / s$,
      precision=2,
      string type,
      column type=Y,
    },
    columns/cores/.style={
      column name=Cores,
      int detect,
    },
    every head row/.style={
      before row={
        \toprule
        & \multicolumn{2}{c}{$Q_1$}
        & \multicolumn{2}{c}{$Q_3$}
        & \multicolumn{2}{c}{$Q_5$}
        & \multicolumn{2}{c}{$Q_7$}\\
        \cmidrule(lr){2-3} \cmidrule(lr){4-5} \cmidrule(lr){6-7} \cmidrule(lr){8-9}
      },
      after row=\midrule,
    },
    every last row/.style={
      after row=\bottomrule
    },
    ]
    {benchmarks-201709/strong-scalability-table.txt}\\[.5em]
  \caption{Run times and efficiencies $\eta$ for strong scalability on IWR compute cluster (416 nodes with 2~x~E5-2630 v3 each, 64 GiB / node, QDR Infiniband)}
  \label{fig:strong-scaling-tab}
\end{table}

\section{Conclusion}
\label{Sec:Conclusion}

We have presented an efficient implementation of a high-order Discontinuous Galerkin discretization for convection-diffusion-reaction
problems that uses sum factorization for optimal algorithmic complexity. Our implementation exploits the inherent structure of this
numerical scheme and delivers more than 50\% of the theoretical peak performance on AVX2-based Intel architectures with fused multiply-add
operations, which make up the majority of compute power in current large-scale compute clusters that is not backed by accelerators. At the
same time, our implementation does not require vectorization over multiple cells or faces and can thus be integrated into existing
frameworks without having to restructure framework code.

For clarity, we have focused in this paper on operator application for scalar problems. In a separate paper, we apply the presented
techniques to a Navier-Stokes solver\cite{piatkowski:2017:navierstokesdg}, while an upcoming paper will investigate efficient preconditioners within this
framework without sacrificing the performance and memory advantages of our approach. At the same time, developing this type of highly
optimized code in C++ is very time consuming and highly affected by the type of problem under investigation and by hardware developments.
For this reason, we are also working on a Python-based code generation tool that can transform a very abstract description of the weak form
into highly optimized code. This tool can already handle complex systems of equations and generate code using the AVX-512 instruction set
for newer versions of Xeon processors as well as the Xeon Phi accelerator architecture. We will describe it in detail in another upcoming paper.

\bibliographystyle{siamplain}
\bibliography{lit.bib}   

\end{document}